\documentstyle[12pt]{article}
\newcommand{\be}{\begin{equation}}
\newcommand{\ee}{\end{equation}}
\newcommand{\bqn}{\begin{eqnarray}}

\newcommand{\eqn}{\end{eqnarray}}

\newcommand{\bd}{\begin{description}}
\newcommand{\ed}{\end{description}}

\newtheorem{stat}{}[section]

\def\bs{\begin{stat}}
\def\es{\end{stat}}
\def\ben{\begin{enumerate}}
\def\een{\end{enumerate}}

\def\bp{\noindent{\bf Proof}  \ \ \ }
\def\ep{\hfill $\Box$}

\makeatletter
\@addtoreset{equation}{section}

\makeatother

\begin{document}

\begin{center}
{\large {\bf ON HAMILTONICITY OF}}
\\[2ex]
{\large {\bf \{CLAW, NET\}-FREE GRAPHS}}
\\[3ex]
{\large {\bf Alexander Kelmans}}
\\[2ex]
{\bf Rutgers University, New Brunswick, New Jersey}
\\[0.5ex]
{\bf University of Puerto Rico, San Juan, Puerto Rico}

\end{center}
 
\begin{abstract}
An {\em $st$-path} is a path with the end-vertices $s$ and $t$. An {\em $s$-path} is a path with an end-vertex $s$.
The results of this paper include necessary and sufficient conditions for a  \{claw, net\}-free graph
$G$ with $s,t \in V(G)$ and $e \in E(G)$ to have
$(1)$
a Hamiltonian  $s$-path,
$(2)$
a Hamiltonian $st$-path,
$(3)$ a Hamiltonian  $s$- and $st$-paths containing $e$
when $G$ has connectivity one, and
$(4)$
a Hamiltonian cycle containing $e$ when $G$ 
is 2-connected.
These results imply that a connected  \{claw, net\}-free graph has a Hamiltonian path and a 2-connected \{claw, net\}-free graph has a Hamiltonian cycle \cite{DGJ}. 
Our proofs of $(1)$-$(4)$ are shorter  than the proofs of their corollaries in \cite{DGJ}, and 
provide
polynomial-time algorithms for solving the corresponding Hamiltonicity problems.
\\[0.5ex]

\end{abstract}

\section{Introduction}

\indent

We consider simple undirected graphs. 
All notions on graphs that are  not defined here can be  found 
in \cite{D,W}.

A graph $G$ is called $H$-{\em free}
if $G$ has no induced subgraph isomorphic to a graph $H$.
A {\em claw} is a graph having exactly four vertices and exactly three edges that are incident to a common vertex.
A claw can be drawn as the letter $Y$.
A {\em net} is a graph obtained from a triangle  by attaching 
to each vertex a new dangling edge.

There are many papers devoted to the study of
Hamiltonicity of claw-free graphs, and, in particular, 
\{claw, net\}-free graphs
(e.g. \cite{BDK,DGJ,FFR,K1,Kcl,K,R,S}).
The maximum independent vertex set problem for 
\{claw, net\}-free graphs was studied in \cite{HMW}.
In this paper we establish some new Hamiltonicity results 
on \{claw, net\}-free graphs.

An {\em $st$-path} is a path with the end-vertices 
$s$ and $t$.
An {\em $s$-path} is a path with an end-vertex $s$.
Let $G$ be a \{claw, net\}-free graph,
$s, t \in V(G)$, $s \ne t$, and $e \in E(G)$.
The results of this paper include necessary and sufficient conditions for $G$ to have:
\\[0.5ex] 
a Hamiltonian $s$-path (see  {\bf \ref{xy-Ham'}} and 
{\bf \ref{2conYNstr}} below),
\\
a Hamiltonian $st$-path when $G$ has connectivity 
one (see {\bf \ref{xy-Ham'}}),
\\
a Hamiltonian  $st$-path containing $e$
if $G$ has connectivity one ({\bf \ref{cn=1SeT}}),
\\
a Hamiltonian  $s$-path containing $e$ when $G$ 
has connectivity one ({\bf \ref{cn=1Se}}), and
\\
a Hamiltonian cycle containing $e$ when $G$ is 
2-connected
({\bf \ref{2conYNstr}}).
\\[0.5ex]
\indent
From the above mentioned results we have the following corollaries.
\bs  {\em \cite{DGJ} (Corollary of {\bf \ref{xy-Ham'}})}
\label{conYN}
$~$Every connected \{claw, net\}-free graph has 
a Hamiltonian path.
\es

\bs  {\em \cite{DGJ} (Corollary of {\bf \ref{2conYNstr}})}
\label{2conYN}
$~$Every 2-connected \{claw, net\}-free graph has 
a Hamiltonian cycle.
\es

Our proofs of {\bf \ref{xy-Ham'}} and {\bf \ref{2conYNstr}}
are shorter and more natural than the proofs of their corollaries {\bf \ref{conYN}} and
{\bf\ref{2conYN}}  in \cite{DGJ}.
They also provide
polynomial time algorithms for solving the corresponding
Hamiltonian problems for \{claw, net\}-free graphs.
In \cite{BDK} a linear time algorithm was given for finding 
a Hamiltonian path and a Hamiltonian cycle (if any exist) in 
a \{claw, net\}-free graph.

The known results on 3-connected \{claw, net\}-free graphs
include the following.
\bs  {\em \cite{S} }
\label{3conHamcon}
A 3-connected \{claw, net\}-free graph 
has a Hamiltonian $xy$-path for every two distinct vertices $x$ and $y$.
\es

\bs  {\em \cite{K} }
\label{3conHam2edges}
Let $G$ be a \{claw, net\}-free graph.
If $G$ is 3-connected, then every two non-adjacent edges 
in $G$ belong to a Hamiltonian cycle.
If $G$ is 4-connected, then every two edges in $G$ belong 
to a Hamiltonian cycle.
\es

\bs  {\em \cite{K} }
\label{3conSeTpath}
Let $G$ be a 3-connected \{claw, net\}-free graph,
$e = uv \in E(G)$, and $s,t \in V(G)$, $s \ne t$.
Then $G$ has a Hamiltonian $st$-path containing $e$ 
if and only if either $\{s,t\} \cap \{u,v\} = \emptyset $
or $\{s,t\} \setminus \{u,v\} = z \in V(G)$
and $G - \{z,u,v\}$ is connected.
\es

\bs  {\em \cite{K} }
\label{LhamconKconYN}
Let $G$ be a $k$-connected \{claw, net\}-free graph, 
$k \ge 3$,
$L_1$ and $L_2$ two disjoint paths in $G$,
$|V(L_1)| + |V(L_2)| \le k$,
and $x_1$, $x_2$ the end-vertices of $L_1$, $L_2$, respectively.
Then the following are equivalent:
\\[0.5ex]
$(c1)$ $G$ has a Hamiltonian $x_1x_2$-path containing $L_1$ and $L_2$,
\\[0.5ex]
$(c2)$ $G$ has a Hamiltonian $z_1z_2$-path containing $L_1$ and $L_2$
for every end-vertices $z_1$, $z_2$ of $L_1$, $L_2$, respectively, and
\\[0.5ex]
$(c3)$ $G - (L_1 \cup L_2)$ is connected.
\es

\bs  {\em \cite{K} }
\label{LhamKconYN}
Let $G$ be a $k$-connected \{claw, net\}-free graph, 
$k \ge 2$, $L$ a path in $G$, and $|V(L)| \le k$.
Then $G$ has a Hamiltonian cycle containing $L$ 
if and only if  $G - L$ is connected.
\es

Obviously both {\bf \ref {3conHamcon}} and 
{\bf \ref {3conHam2edges}} follow immediately from 
{\bf \ref{3conSeTpath}}.
More results on Hamiltonicity of $k$-connected 
\{claw, net\}-free graphs can be found in \cite{K}.

The results of this paper form a part of a broader picture on Hamiltonicity of  \{claw, net\}-free graphs and were presented at the Discrete Mathematics Seminar at the University of Puerto Rico in November 1999  (see also \cite{K,K1}).

\section{Main notions and notation}
\label{notation}
\indent

We consider undirected graphs with no loops and no parallel edges. We use the following notation:
$V(G)$ and $E(G)$ are the sets of vertices and edges of 
a graph $G$, respectively,
$v(G) = |V(G)|$ and $e(G) = |E(G)|$,
$AvB$ is the union of two graphs $A$ and $B$ having exactly one vertex $v$ in common, and $AvB = Avu$
if $B$ is an edge $vu$.

An $st$-{\em path} ($s$-{\em path})
is a path with the end-vertices $s$ and $t$ (an end-vertex $s$, respectively).
If $a$ and $b$ are vertices of $P$, then $aPb$ denote  
the subpath of $P$ with the end-vertices $a$ and $b$.
A path (a cycle) of $G$ is called {\em Hamiltonian} if it contains each vertex of $G$.
A Hamiltonian path of $G$ is also called a {\em trace} of $G$.
We introduce the term {\em track} of $G$ for a Hamiltonian cycle of $G$.

Let $\kappa(G)$ denote the {\em vertex connectivity} of 
a graph $G$. A  graph $G$ is called $k$-{\em connected} 
if $\kappa(G) \ge k$.

Let $H$ be a subgraph of $G$. We write simply $G - H$ instead of $G - V(H)$. A vertex $x$ of $H$ is called
an {\em inner vertex of}
$H$ if $x$ is adjacent to no vertices in $G - H$, and
a {\em boundary vertex of} $H$, otherwise.
An edge $e$ of $H$ is called an {\em inner edge of} $H$ if
$e$ is incident to an inner vertex of $H$.

A {\em block} of $G$ is either an isolated vertex or a maximal connected subgraph $H$ of $G$ such that $H - v$ is connected for every $v \in V(H)$. 
A block $B$ of $G$ is called an {\em end-block of} $G$ if $B$ has exactly one
boundary vertex, and an {\em inner block}, otherwise.

\section{ The key lemma}
\label{keylemma}
\indent

First we observe the following.
\bs
\label{3ends}
Let $G$ be a graph. The following are equivalent:
\\[0.5ex]
$(a1)$ $G$ has no induced subgraph isomorphic to 
a claw or a net and
\\[0.5ex]
$(a2)$ $G$ has no connected induced subgraph with 
at least three end-blocks.
\es

\bp
Obviously $(a2)\Rightarrow (a1)$.
We prove $(a1)\Rightarrow (a2)$.
If $G$ is \{claw, net\}-free, then $G - x$ is also 
\{claw, net\}-free for every $x \in V(G)$. 
Clearly our claim is true if $v(G) = 1$.
Let $F$ be a counterexample with the minimum number of vertices. 
Then 
$(1)$ every end-block  has exactly one edge, 
$(2)$ $F$ has exactly three end-blocks, 
$(3)$ if $x \in V(F)$ and $F - x$ is connected, then $x$ is 
a leaf, and
$(4)$ $F$ is not a claw and not a net.
By $(2)$ and $(3)$, $F$ is a tree or has exactly one cycle which is a triangle. 
In both cases by $(4)$,
$F$ has a leaf $z$ such that 
$F - z$ is  a smaller counterexample, a contradiction.
\ep
\\[1ex]
\indent
The following lemma is 
useful for analyzing Hamiltonicity of  \{claw, net\}-free graphs.
\bs
\label{L}
Let $G$ be a  \{claw, net\}-free graph and $z \in V(G)$. 
Suppose that $G - z$ has an $xy$-trace $P$ and 
there exists  $e_z = zp \in E(G)$, and so $G$ is connected and $p \in V(P)$.
Let $e_x$ and $e_y$ be the end-edges of $P$.
Then $G$ has an $ab$-trace $Q$ 
such that $\{a,b\} \subset \{x,y,z\}$,
$e_z \in  E(Q)$ and $\{e_x,e_y\} \cap E(Q) \ne \emptyset $.
\es
\noindent
\bp (uses {\bf \ref{3ends}}). 
We  define below a notion of  
a {\em good path} which is a special subpath of path $P$.
Our goal is to show that if $G$ has no required trace, 
then $G$ has a good path and a maximal good path 
is  a subpath of a longer good path in $G$, which is  
a contradiction. 

By the assumption of our claim, $p \in V(P)$.
Let $X = pPx = x_0x_1 \cdots x_{k-1}x_k$ and 
$Y = pPy = y_0y_1 \cdots y_t$, where $x_k = x$, $y_t = y$, 
and $x_0 = y_0 = p$.
Let $M_{r,s} = x_rPy_s$,
$\dot{M}_{r,s}$ denote the subgraph of $G$ induced by $V(M_{r,s})$, and
$\bar{M}_{r,s} = \dot{M}_{r,s} \cup \{x_rx_{r+1}, y_sy_{s+1},zp\}$.
\\[0.8ex]
\indent
A subpath $M_{r,s}$ is called {\em good} if
\\[0.5ex]
{\bf (x1)} $\dot{M}_{r,s}$ has a $py_s$-trace containing 
$x_{r-1}x_r$,
\\[0.5ex]
{\bf (y1)} $\dot{M}_{r,s}$ has a $px_r$-trace containing $y_{s-1}y_s$,
\\[0.5ex]
{\bf (x2)} if $x_r \ne x$, then for every  $v \in V(M_{r,s}) \setminus x_r$, the graph 
$\dot{M}_{r,s} \cup  \{x_rx_{r+1},x_{r+1}v\}$  
obtained from $\dot{M}_{r,s}$ by
adding the edge $x_rx_{r+1}$ and a new edge $x_{r+1}v$ 
has a $py_s$-trace containing the path $x_rx_{r+1}v$,
\\[0.5ex]
{\bf (y2)} if $y_s \ne y$, then for every
 $v \in V(M_{r,s}) \setminus y_s$, the graph
$\dot{M}_{r,s} \cup \{y_sy_{s+1},y_{s+1}v\}$ obtained from $\dot{M}_{r,s}$ by adding
the edge $y_sy_{s+1}$ and a new edge $y_{s+1}v$ 
 has a $px_r$-trace containing the path $y_sy_{s+1}v$, 
 and
\\[0.5ex]
{\bf (z)} for every $v \in V(M_{r,s}) \setminus p$, the graph 
$\dot{M}_{r,s} \cup \{zp,zv\}$
obtained from $\dot{M}_{r,s}$ by adding the edge $zp$ 
and a new edge $zv$ has an $x_ry_s$-trace
(which clearly contains $e_z = zp$ and $zv$).
\\[0.8ex]
\indent
If $p \in \{x,y\}$ or $\{x_1z,y_1z \} \cap E(G) \ne \emptyset$,
then clearly $G$ has a required trace. 
Therefore let $p \not \in \{x,y\}$ and $\{x_1z,y_1z \} \cap E(G) =  \emptyset $. Since $G$ has no induced claws, 
the claw in $G$ with the edge set $\{px_1,py_1,pz\}$ is 
not induced, and therefore $x_1y_1 \in E(G)$.

Clearly $\dot{M}_{1,1}$ is a triangle and 
$V(\dot{M}_{1,1}) = \{p,x_1,x_2\}$. 
Now it is easy to check that $M_{1,1}$ is a good path.
Let $M_{r,s}$ be a maximal good path.
Put $A = \{e_x, e_y, e_z\}$.
\\[0.5ex]
${\bf (p1)}$
Suppose  that $x_r = x$. By {\bf (x1)}, $\dot{M}_{r,s}$ has 
a $py_s$-trace $L$ containing $x_{r-1}x_r$.
Then $zpLy_sPy$ is a $yz$-trace in $G$ containing $A$.
Similarly, if $y_s = y$, then  $G$ has an $xz$-trace  containing $A$. 
\\[0.5ex]
${\bf (p2)}$ 
Now suppose that $x_r \ne x$ and $y_s \ne y$.
Then the subgraph $\bar{M}_{r,s}$ of $G$ has at least three end-blocks.
Since $G$ is  \{claw,~net\}-free, by {\bf \ref{3ends}}, there exists an edge $ab$ in $G$ such that
$a \in \{x_{r+1},y_{s+1},z\}$ and $b \in V(\bar{M}_{r,s} - a)$.
\\[0.5ex]
${\bf (p2.1)}$ Suppose that $a = z$ and $b \in V(M_{r,s})$. 
By {\bf {\bf (z)}},
$\bar{M}_{r,s} \cup zb$ has an $x_{r+1}y_{s+1}$-trace $L$
containing $e_z$.  Then $xPx_{r+1}Ly_{s+1}Py$ is an 
$xy$-trace in $G$ containing $A$.
\\[0.5ex]
${\bf (p2.2)}$ Suppose that $a = z$ and 
$b \in \{x_{r+1},y_{s+1}\}$. 
By symmetry, we can assume that $b = x_{r+1}$.
By {\bf (x1)}, $\dot{M}_{r,s}$ has a $py_s$-trace $L$.
Then $P' = xPx_{r+1}zpLy_sPy$ is an $xy$-trace in $G$.
If $x \ne x_{r+1}$, then $P'$ contains $A$. If  $x = x_{r+1}$, 
then $P'$ contains $A \setminus e_x$.
\\[0.5ex]
${\bf (p2.3)}$ Now suppose that $a \in \{x_{r+1},y_{s+1}\}$ 
and $b \ne z$.
By symmetry, we can assume that  $a = x_{r+1}$.
Then $b \in  V(M_{r,s} - x_r) \cup y_{s+1}$.
\\[0.5ex]
${\bf (p2.3.1)}$ Suppose that $x_{r+1} = x$.

Suppose that $b \ne y_{s+1}$.
By {\bf (x2)}, $M_{r,s} \cup xb$ has a $zy_s$-trace $L$
containing $e_x = x_rx_{r+1}$.
Then $zpLy_sy_{s+1}Py$ is a $yz$-trace in $G$ containing $A$.

Now suppose that $b = y_{s+1}$.
By {\bf (y1)}, $\dot{M}_{r,s}$ has a $\{p,x_r\}$-trace $L$.
Then $P' = zpLx_rx_{r+1}y_{s+1}Py$ is a $zy$-trace in $G$.
If $y_{s+1} \ne y$, then $P'$ contains $A$.
If $y_{s+1} = y$, then $P'$ contains $A - e_y$.
\\[0.5ex]
${\bf (p2.3.2)}$ Now suppose that $x_{r+1} \ne x$. 
Our goal is to show that
\\[0.5ex]
(c1) if $b \ne y_{s+1}$, then $M' = M_{r+1,s}$ is a good  path and
\\[0.5ex]
(c2) if $b = y_{s+1}$ (i.e. $x_{r+1}y_{s+1} \in E(G)$), then
$M' = M_{r+1,s+1}$ is a good  path.
\\[0.5ex]
\indent This will lead to  a contradiction because  $M_{r,s} \subset M'$, and
therefore a good path  $M_{r,s}$ will not be maximal.
We recall that we consider the case when $x_r \ne x$  and 
$y_s \ne y$.
\\[0.5ex]
\indent
{\sc Case} $(c1)$.
Suppose that $b \ne y_{s+1}$.
We 
want to 
prove that $M_{r+1,s}$ is a good path.
\\[0.5ex]
${\bf (p.x1)}$ Let us show that $M_{r+1,s}$ satisfies 
{\bf (x1)}.
By {\bf (x2)} for $M_{r,s}$, the graph
$\dot{M}_{r,s} \cup \{x_rx_{r+1},x_{r+1}b\}$ has  
a $py_s$-trace $L$ containing the path $x_rx_{r+1}b$. 
Then $L$ is also a $py_s$-trace in  $\dot{M}_{r+1,s}$ containing $x_rx_{r+1}$.
\\[0.5ex]
${\bf (p.y1)}$
Let us show that $M_{r+1,s}$ satisfies {\bf (y1)}.
By {\bf (y1)} for $M_{r,s}$, the graph $\dot{M}_{r,s}$ has 
a $px_r$-trace $L$ containing $y_{s-1}y_s$. 
Then $pLx_rx_{r+1}$ is a $px_{r+1}$-trace in $\dot{M}_{r+1,s}$ containing $y_{s-1}y_s$.
\\[0.5ex]
${\bf(p.x2)}$ Let us show that $M_{r+1,s}$ satisfies {\bf (x2)}.

Consider graph $Q_v = \dot{M} \cup \{x_{r+1}x_{r+2},   x_{r+2}v\}$, where $v \in V(M_{r+1,s}) \setminus x_{r+1}$.

Suppose that $v \ne x_r$. By {\bf (x2)} for $M_{r,s}$, graph
$\dot{M}_{r,s} \cup  \{x_rx_{r+1}, vx_{r+1}\}$ has 
a $py_s$-trace $L$ containing the path $x_rx_{r+1}v$. 
Then $(L - vx_{r+1}) \cup (x_{r+1}x_{r+2}v)$ is a 
$py_s$-trace in $Q_v$ containing path $x_{r+1}x_{r+2}v$.

Now suppose that $v = x_r$. By ${\bf (p.x1)}$, $M_{r+1,s}$ satisfies
{\bf (x1)}, i.e. graph $\dot{M}_{r+1,s}$ has 
a $py_s$-trace $L$ containing $x_rx_{r+1}$. 
Then $(L - x_rx_{r+1}) \cup (x_{r+1}x_{r+2}x_r)$ is  
a $py_s$-trace containing path $x_{r+1}x_{r+2}v$.
\\[0.5ex]
${\bf(p.y2)}$ Let us show that $M_{r+1,s}$ satisfies {\bf (y2)}.

Consider graph 
$Q_v = \dot{M}_{r+1,s} \cup  \{y_sy_{s+1}, vy_{s+1}\}$, 
where $v \in V(M_{r+1,s}) \setminus y_s$.
By {\bf {\bf (y2)}} for $M_{r,s}$, graph
$\dot{M}_{r,s} \cup \{y_sy_{s+1}, vy_{s+1}\}$ has 
a $px_r$-trace $L$ containing path $y_sy_{s+1}v$.
Then  $x_{r+1}x_rLz$ is a $\{p,x_{r+1}\}$-trace in $Q_v$ containing path $y_sy_{s+1}v$.
\\[0.5ex]
${\bf (p.z)}$ Let us show that $M_{r+1,s}$ satisfies {\bf (z)}.

Consider graph $Q_v = M_{r+1,s} \cup \{zp,zv\}$,
where $v \in V(M_{r+1,s}) \setminus p$.

Suppose that $v \in V(M_{r,s}) \setminus p$. By {\bf (z)} for $M_{r,s}$, graph $M_{rs} \cup \{zp,zv\}$ has 
an $x_ry_s$-trace $L$. 
Then $x_{r+1}x_rLy_s$ is an $x_{r+1}y_s$-trace
in $M_{r+1,s} \cup \{zp,zv\}$.

Now suppose that $v = x_{r+1}$.
By {\bf (x1)} for $M_{r,s}$, graph
$\dot{M}_{r,s}$ has a $py_s$-trace $L$.
Then
$x_{r+1}zpLy_s$ is an $x_{r+1}y_s$-trace in $Q_v$.
\\[0.5ex]
\indent
{\sc Case} $(c2)$.
Now suppose that $b = y_{s+1}$. 
We want to prove that $M_{r+1,s+1}$ is a good path. 
By symmetry, it suffices to proof  
that $M_{r+1,s+1}$ satisfies 
{\bf (x1)},  {\bf (x2)}, and  {\bf (z)}. Let us proof  {\bf (x1)}.
By {\bf (y1)} for $M_{r,s}$, graph
$\dot{M}_{r,s}$ has a  $px_r$-trace $L$.
Then $pLx_rx_{r+1}y_{s+1}$ is a $py_{s+1}$-trace in $\dot{M}_{r+1,s+1}$ containing $x_rx_{r+1}$.
The proof of {\bf (x2)} and  {\bf (z)} is similar to  
{\sc Case} $(c1)$.
\ep

\section{More on \{claw, net\}-free graph Hamiltonicity}
\label{Hamilton}

\indent

Lemma {\bf \ref{L}} allows to give an easy proof of the following strengthening
of {\bf \ref{conYN}}.
\bs
\label{conYNstr1} Let $G$ be a connected \{claw, net\}-free graph.
Then
\\[0.5ex]
$(a1)$ $G$ has a trace and
\\[0.5ex]
$(a2)$ if $sz \in E(G)$ and $G - z$ is connected, then
$sz$ belongs to a trace of $G$.
\es

\bp (uses {\bf \ref{L}}).
We prove our claim by induction on $v(G)$.
The claim holds if $v(G) = 1$.
Since $G$ is connected, there exists $z \in V(G)$
such that $G - z$ is also connected.
Let $sz \in E(G)$.
Since $G$ is \{claw, net\}-free, clearly $G - z$ is also
\{claw, net\}-free.
Therefore by the induction hypothesis, $G - z$ has a trace.
Then by {\bf \ref{L}}, $G$ has a trace containing $sz$.
\ep
\\[1ex]
\indent
Here is another strengthening of {\bf \ref{conYN}} for graphs of connectivity one.
\bs
\label{xy-Ham}
Let $G$ be a connected \{claw, net\}-free graph, 
$G = AaHbB$, where $A$ and $B$ are
end-blocks of $G$. Let $a' \in V(A - a)$, $b' \in V(B - b)$, 
and $a'x$ be an edge of $A$ such that if $v(A) \ge 3$,
then $x$ is an inner vertex of an end-block of $G - a'$.
Then
\\[0.5ex]
$(a1)$ there exists an $a'b'$-trace in $G$ and, moreover,
\\[0.5ex]
$(a2)$ there exists an $a'b'$-trace in $G$ containing edge $a'x$.
\es

\bp 
We prove our claim by induction on $v(G)$.
If $v(G) = 3$, then our claim is obviously true.
\\[0.5ex]
{\bf (p1)} Suppose that $v(A) \ge 3$.
Then $A$ is 2-connected.
Let $A' = A - a'$ and $G' = G - a'$.
Then $G' = A'aHbB$ and $G'$ is connected.
Since $G$ is \{claw, net\}-free, $G'$ is also
\{claw,~net\}-free.
Since $v(G') < v(G)$, by the induction hypothesis,
$G'$ has an $xb'$-trace $P$. Then 
$a'xPb'$ is an $a'b'$-trace in $G$
containing $a'x$.
\\[0.5ex]
{\bf (p2)} Now suppose that $v(A) = 2$.
Then $a'x = a'a$ and there is $b'z \in E(B)$ such that $z$ 
is an inner vertex of an end-block in $G - b'$.
Hence by the arguments, similar to those in {\bf (p1)}, 
$G$ has an $a'b'$-trace in $G$
containing $a'x$ (as well as $b'z$).
\ep
\\[1ex]
\indent
From {\bf \ref{xy-Ham}} we have, in particular:
\bs
\label{xy-Ham'}
Let $G$ be a \{claw, net\}-free graph, $v(G) \ge 3$, $\kappa(G) = 1$, and $s, t \in V(G)$.
Then $G$ has an $st$-trace if and only if
$s$ and $t$ are inner vertices of different end-blocks of $G$.
\es

From {\bf \ref{conYNstr1}} and {\bf \ref{xy-Ham}}
it is easy to obtain the following stronger result.

\bs
\label{conYNstr2}
Let $G$ be a connected \{claw, net\}-free graph having 
$k \ge 2$ blocks.
Let $A_j$, $j \in \{1,2\}$,  be an end-block of $G$, 
$a'_j$ the boundary vertex
of $A_j$, $a_j \in A_j - a'_j$, and $\alpha _j \in E(A_j)$.
Let $B_i$ be an inner block of $G$ and  
$\beta _i \in E(B_i)$.
Let $U = \{\alpha _1, \alpha _2\} \cup \{\beta _i: i = 1, \ldots , k-2\}$.
Suppose that
\\[0.5ex]
$(h1)$
$\alpha _j = a_jx_j$ is such that
if $v(A) \ge 3$, then $x_j$ is an inner vertex of an end-block
of $A_j - a'_j$, $j \in \{1,2\}$, and
\\[0.5ex]
$(h2)$
$\beta _i$ is an inner edge of $B_i$, if $v(B_i) \ge 3$,
$i \in \{ 1, \ldots , k-2\}$.

Then $G$ has an $a_1a_2$-trace containing $U$.
\es

\bp (uses {\bf \ref{conYNstr1}} and {\bf \ref{xy-Ham}}).
Since $G$ is connected, for every end-block $A_j$ of $G$ there is an edge
$a'_jp_j \in E(G) \setminus E(A_j)$.
Similarly, for every inner block $B_i$ of $G$ there are edges
$b_iq_j, b'_iq'_j \in E(G) \setminus E(B_i)$, where
$b_i$  and $b'_i$ are the boundary vertices of $B_i$.
Let $\bar{A}_j = A_ja'_jp_j$ and $\bar{B}_i = q_ib_iB_ib'_iq'_j$.
Then all $\bar{A}_j$'s and $\bar{B}_i$'s are induced subgraphs of $G$ and,
therefore, are \{claw, net\}-free.
By {\bf \ref{conYNstr1}}, each $\bar{B}_i$ has
a trace $q_ib_iQ_ib'_iq'_j$ containing $\beta _i$.
By {\bf \ref{xy-Ham}}, each $\bar{A}_j$ has
a trace $a_jP_ja'_jp_j$ containing $\alpha _j$.
Then $P_1 \cup Q_1 \ldots Q_{k-2}\cup P_2$ is
an $a_1a_2$-trace containing $U$.
\ep
\\[1ex]
\indent
Let ${\cal L}$ denote the set of 4-tuples $(G, s,t, uv)$ such that $G$ is a graph, $\{s, t\} \subseteq  V(G)$, $s \ne t$, $uv \in E(G)$, and either $(1)$ $\{s,t\}$ does not meet one of the components of $G - \{u,v\}$ or $(2)$ $\{s,t\} \cap \{u,v\} \ne \emptyset $, say $t = u$, and
either $G - \{s,v\}$ is not connected and the component containing $t$ has at least two vertices or
there is $x \in V(G - \{u,v\})$ such that $\{s,v\}$ avoids one of the components of $G - \{t, x\}$.

Obviously, if $G$ has an $st$-trace containing $uv$, 
then $(G, s,t, uv) \not \in {\cal L}$.
We will see that for \{claw, net\}-free graphs of connectivity one the converse is also true.

\bs
\label{e-Ham,Gy}
Let $G$ be a connected graph, $s \in V(G)$, and $xsG$ 
be a  \{claw, net\}-free graph.
Let $C$ be the end-block of $vsG$ distinct from $xs$, 
$c$ the boundary vertex of $C$,
$t \in V(C - c)$, and $uv \in E(G)$.
Then $G$ has an $st$-trace containing  
$uv$ if and only if $(G, s,t, uv) \not \in {\cal L}$.
\es

\bp (uses {{\bf \ref{xy-Ham}} and {\bf \ref{conYNstr2}}).
By the above remark, it is sufficient to show that
$(G,s,t,uv) \not \in {\cal L}$ implies that $G$ has 
an $st$-trace containing  $uv$.
We prove our claim by induction on $v(G)$.
If $uv \not \in E(C)$ or $V(C) = \{u,v\}$, then our claim follows from {\bf \ref{conYNstr2}}. Therefore let $uv  \in E(C)$.
In particular, if $v(C) = 2$, then our claim is true.
Therefore let $v(C) \ge 3$, and so $C$ is 2-connected.
Let $G' = G - t$ and $C' = C - t$, and so $C'$ is connected.
\\[0.5ex]
${\bf (p1)}$ Suppose that $G - \{u,v\}$ is not connected.
Since $(G, s,t, uv) \not \in {\cal L}$, vertices $s$ and $t$ belong in $G - \{u,v\}$
to different components, say $S$ and $T$, respectively.
Since $C$ is 2-connected, $\bar{T} = T \cup uv$ 
is also 2-connected.
\\[0.5ex]
${\bf (p1.1)}$ Suppose that $v(T) = 1$, i.e. $V(T) = \{t\}$.
Then $tu$ is an end-block of $G - v$.
Since $xsG$ is \{claw, net\}-free, by {\bf \ref{xy-Ham}},
$G - v$ has an $st$-trace $sPut$.
Then $sPuvt$ is an $st$-trace in $G$ containing $uv$.
\\[0.5ex]
${\bf (p1.2)}$ Now suppose that $v(T) \ge 2$.  
Since $\bar{T}$ is 2-connected, either $\bar{T} - t$ is 
2-connected or $t$ is adjacent in $G$ to an inner vertex $z$ of the end-block of $\bar{T} - t$ avoiding $uv$.
In both cases, $(G', s,z, uv) \not \in {\cal L}$, and so by the induction hypothesis, $G'$  has a $sz$-trace $P$ containing $uv$.
Then $sPzt$ is an $st$-trace containing $uv$.
\\[0.5ex]
${\bf (p2)}$ Now suppose that $G - \{u,v\}$ is connected.
Since $(G, s,t, uv) \not \in {\cal L}$,  $\{u,v\} \ne \{s,t\}$. 
Since $C$ is 2-connected, $t$ is adjacent to an inner vertex $z$ of the end-block $B$ of $xsG'$ which avoids $x$.
If $t \in \{u,v\}$, say $t = a$, then since 
$(G, s,t, uv) \not \in {\cal L}$, $v$ is an inner vertex of $B$.
Then by {\bf \ref{xy-Ham}}, $G'$ has an $sv$-trace $P$, 
and so $sPba$ is an $st$-trace containing $uv$.
So let $t \not \in \{u,v\}$.
Let $D$ be the block of $G'$ containing $uv$.
If $D \ne B$, then since $(G, s,t, uv) \not \in {\cal L}$, 
also $(G', s,z, uv)  \not \in {\cal L}$, and so 
by the induction hypothesis, $G'$  has a $sz$-trace $P$ containing $uv$.
If $D = B$, then $(G, s, z, uv) \not \in {\cal L}$ 
because $G$ has no induced claw centered at $z$. 
So again by the induction hypothesis, $G'$  has 
a $sz$-trace $P$ containing $uv$.
In both cases $sPzt$ is an $st$-trace in $G$ containing $uv$.
\ep
\\[1ex]
\indent
From {\bf \ref{conYNstr2}} and {\bf \ref{e-Ham,Gy}} we have:
\bs
\label{cn=1SeT}
Let $G$ be a \{claw, net\}-free graph, 
$v(G) \ge 3$, 
$\kappa(G) = 1$, $e \in E(G)$, 
and $\{s,t\} \in V(G)$, $s \ne t$.
Then $G$ has an $st$-trace containing $e$ 
if and only if
$s$ and $t$ are inner vertices of different end-blocks of 
$G$ and $(G, s,t, e) \not  \in {\cal L}$.
\es

From {\bf \ref{cn=1SeT}} we have:
\bs
\label{cn=1Se}
Let $G$ be a \{claw, net\}-free graph, $v(G) \ge 3$, 
$\kappa(G) = 1$, $s \in V(G)$, and $e \in E(G)$.
Then $G$ has an $s$-trace containing $e$ if and only if
$s$ is an inner vertex of an end-block in $G$ and
$(G, b,s, e) \not  \in {\cal L}$, where $b$ is the boundary vertex of the end-block  avoiding $s$.
\es

From {\bf \ref{conYNstr2}} and {\bf \ref{cn=1SeT}} we have
the following strengthening of {\bf \ref{conYNstr2}}.

\bs
\label{conYNstr3}
Let $G$ be a connected \{claw, net\}-free graph having 
$k \ge 2$ blocks.
Let $A_j$, $j \in \{ 1,2\}$,  be an end-block of $G$, $a'_j$ 
the boundary vertex
of $A_j$, $a_j \in A_j - a'_j$, and $\alpha _j \in E(A_j)$.
Let $B_i$ be an inner block of $G$ and  $\beta _i \in E(B_i)$.
Let 
$U = \{\alpha _1, \alpha _2\} \cup \{\beta _i: i = 1, \ldots , 
k-2\}$.
Then $G$ has an $a_1a_2$-trace containing $U$ 
if and only if
\\[0.5ex]
$(c1)$ 
$(A_j, a_j, a'_j, \alpha _j) \not \in {\cal L}$, 
$j \in \{1,2\}$ and
\\[0.5ex]
$(c2)$ $\beta _i$ is an inner edge of $B_i$ if $v(B_i) \ge 3$,
$i \in \{1, \ldots , k-2\}$.
\es

Let ${\cal E}$ denote the set of tuples $(G, e)$ such that
$G$ is a 2-connected graph, $e = x_1x_2 \in E(G)$, 
$G = x_1G_1x_2G_2x_1$, and $G_i \cup x_1x_2$ is 
3-connected or a triangle for some $i \in \{1,2\}$.
 
Obviously, if $e$ belongs to a track of $G$, 
then $(G,e)\not \in {\cal E}$. 
The following strengthening of {\bf \ref{2conYN}} shows, 
in particular, that for 2-connected \{claw, net\}-free graphs the converse  is also true.
\bs
\label{2conYNstr}
Let $G$ be a 2-connected \{claw, net\}-free graph and 
$e = pz \in E(G)$.
Then
\\[0.5ex]
$(a1)$ $G$ has a track,
\\[0.5ex]
$(a2)$ the following are equivalent:
\\[0.5ex]
\indent
$(c1)$ $e$ belongs to a track of $G$,
\\[0.3ex]
\indent
$(c2)$ $(G,e)\not  \in {\cal E}$, and
\\[0.5ex]
$(a3)$ if $(G, e) \in {\cal E}$, then for every inner vertices $s$, $t$ of the two different blocks $S$ and $T$ 
of $G - z$ that contain $p$, there is an 
$st$-trace of $G$ containing $e$.
\es
\bp (uses {\bf \ref{L}} and  {\bf \ref{xy-Ham}} $(a1)$).
As we mentioned above, $(c1)\Rightarrow (c2)$.
\\[0.5ex]
${\bf (p1)}$
We prove $(a1)$ and $(c2) \Rightarrow (c1)$ by induction on $v(G)$.
The claim holds, if $v(G) = 3$ or
$G$ is a cycle. Therefore let $v(G) \ge 4$ and $G$  not 
a cycle. By $(c2)$, $(G,pz) \not  \in {\cal E}$.
\\[0.5ex]
${\bf (p1.1)}$ Suppose that  $G - z$ is 2-connected.
Since $G$ is \{claw, net\}-free, clearly $G - z$ is also
\{claw, net\}-free. Therefore by the induction hypothesis, 
$G - z$ has a track $C$, and so $p \in V(C)$.
Since $G$ is 2-connected, there is a vertex $c$
in $C$ distinct from $p$ and adjacent to $z$.
Let $x$ and $y$ be the two vertices adjacent to $c$ in $C$.
Then $G' = G - c$ satisfies the assumptions of 
{\bf \ref{L}}, namely,
$G'$ is connected and $P = C - c$ is an $xy$-trace of 
$G' - z$.
By {\bf \ref{L}}, $G'$ has an $st$-trace $L$ such that 
$e \in E(L)$ and $\{s,t\} \subset\{x,y,z\}$.
Since $c$ is adjacent to $x$, $y$, and $z$, clearly $csLtc$ 
is a track of $G$ containing $e$.
\\[0.5ex]
${\bf (p1.2)}$
Now suppose that $G - z$ is not 2-connected.
Let $G - z = AaHbB$, where $A$ and $B$ are 
end-blocks of $G$.
Since $(G,pz) \not  \in {\cal E}$, $p$ is an inner vertex of an end-block, say $p \in V(A - a)$.
Since $G$ is 2-connected, $(G,qz) \not  \in {\cal E}$
for some $q \in V(B - b)$.
By {\bf \ref{xy-Ham}} $(a1)$, $G - z$ has a $pq$-trace $P$. 
Then $zpPqz$ is a track in $G$ containing  $e = pz$.
\\[0.5ex]
${\bf (p2)}$ Now we prove $(a3)$.  
Let $(G,pz) \not  \in {\cal E}$. Then
$G - z = SpTbB$, where $S$ is an end-block and $T$ is 
a block of  $G - z$. Let $s$ and  $t$ be inner vertices of $S$ and $T$, respectively. Since $G$ is 2-connected,
$G - S$ is connected. Since $G$ is claw-free, $T - S$ is an end-block of $G - S$, and so $t$ and $z$ are inner vertices of different 
end-blocks of $G - S$.
By {\bf \ref{xy-Ham}} $(a1)$, $S$ has an 
$sp$-path $P$ and $G - S$ has a  
$zt$-trace $Q$. Then  $sPpzQt$ is an $st$-trace of $G$ containing $e$.
\ep
\\[1ex]
\indent
From {\bf \ref{2conYNstr}} we have, in particular:
\bs
\label{2-con,einHampath}
Let $G$ be a 2-connected \{claw, net\}-free graph.
Then every edge in $G$ belongs to a trace of $G$.
\es

In \cite{K4} we gave a structural
characterization of so-called `closed' \{claw, net\}-free graphs.
This structure theorem together with the known properties of the Ryj\'{a}\u{c}ek closure \cite{R} can be used to provide alternative proofs for  some of the above Hamiltonicity results.
In \cite{Kcl} we describe some graph closures that are stronger than the closure in \cite{R}
and that can be applied to graphs having some induced claws. 
These results can be used to extend the picture, described in this paper, for a wider class of graphs.

\end{document}